\documentclass[12pt,letter]{amsart}
\usepackage[latin1]{inputenc}
\usepackage{amsmath,amsthm,amssymb,latexsym}
\usepackage[T1]{fontenc}
\usepackage{libertine} 
\usepackage{euler} 
\usepackage{eucal} 
\usepackage{color}
\usepackage{ulem}
\usepackage{soul}
\usepackage{graphicx}
\usepackage{pb-diagram}
\usepackage{pst-3d,pst-coil,pst-eps,pst-fill,pst-grad,pst-math,pst-text,pst-tree}

\newcommand\C[1]{\mathcal{#1}}

\newtheorem{proposition}{Proposition}[section]
\newtheorem{theorem}[proposition]{Theorem}

\newtheorem{lemma}[proposition]{Lemma}

\newtheorem{fact}[proposition]{Fact}

\theoremstyle{definition}
\newtheorem{definition}[proposition]{Definition}

\newtheorem{remark}[proposition]{Remark}

\def\monster{\mathbb{M}}
\def\E{\varepsilon}
\def\dist{\mathbf{d}}

\def\ordencea{\prec_{\C{K}}}
\def\ordenceac{\succ_{\C{K}}}
\newcommand{\eop}[1]{
\hspace{10mm} \vspace{-6mm}
\begin{flushright}
\qedsymbol$_{\text{#1}}$\\ \ \\
\end{flushright}
}
\newenvironment{prueba}[1][{\it Proof}]{\noindent {\it #1.} }{}
\newcommand{\bdem}[1][Proof]{\begin{prueba}[#1]}
\newcommand{\edem}[1][]{\eop{#1}
\end{prueba}}
\def\bsdem{\begin{prueba}[Reference]}\def\bsindem{\begin{proof}[\ ]}

\def\tuple{\overline}
\def\rest{\upharpoonright}

\mathchardef\mhyphen="2D
\newcommand{\gatp}{{\sf ga\mhyphen tp}}
\newcommand{\gaS}{{\sf ga\mhyphen S}}

\newcommand{\sea}{\mathfrak{C}}

\begin{document}
\title[Around the set-theoretical consistency of d-tameness $\cdots$]{Around the set-theoretical consistency of d-tameness of Metric
  Abstract Elementary Classes.}
\author[W. Boney]{Will Boney}
\author[P. Zambrano]{Pedro Zambrano}
\thanks{W. Boney: Mathematics Department Harvard University, Cambridge, MA. Email: wboney@math,harvard.edu\\
P. Zambrano: Departamento de Matem\'aticas Universidad Nacional de Colombia, AK 30 $\#$ 45-03 c\'odigo postal 111321, Bogota, Colombia. Email: phzambranor@unal.edu.co\\
This material is based upon work while the first author was supported by the National
Science Foundation under Grant No. DMS-1402191.}

\date{\today}
\begin{abstract}
In this paper, we prove that if $\kappa$ is a almost strongly compact cardinal, then any MAEC with L\"owenheim-Skolem number below $\kappa$ is $<\kappa$-d-tame, using ideas from  \cite{Bo14}
\end{abstract}
\maketitle

\section{Motivation}

S. Shelah  and J. Stern proved in \cite{ShSt78} that the Hanf number of the first order theory of Banach spaces is very high, in fact it has the behavior of the second order logic of binary relations. In the 60's and 70's, C.C. Chang and J. Keisler, and independently W. Henson, began the development of a suitable logic to understand Banach spaces with similar properties as first order logic. More recently, this approach was rediscovered (in the 90's and in the 2000's) by W. Henson, J. Iovino, I. Ben-Yaacov, A. Usvyatsov et al, providing a more workable framework (continuous first-order logic, see \cite{CoMon}) to understand classes of complete metric structures, where all logical symbols are interpreted as uniformly continuous functions in {\it complete} metric spaces. 

Metric Abstract Elementary Classes (shortly, MAECs) corresponds to a suitable abstract model-theoretical notion due to \r{A}. Hirvonen and T. Hyttinen (\cite{HiHy}) for studying non-elementary (in the sense of continuous logic) classes of complete metric spaces, in a similar way as Abstract Elementary Classes (AECs). The importance of tame AECs lies in the fact that R.  Grossberg and M. VanDieren provided in \cite{GrVa} a partial answer of the Shelah's categoricity transfer conjecture in this setting. W. Boney proved in \cite{Bo14} that under the existence of a proper class of strongly compact cardinals all AECs are tame, providing a proof of the set-theoretical consistency of the Shelah's categoricity transfer conjecture. P. Zambrano studied in \cite{Za12} a metric version of tameness (d-tameness) in MAECs, proving a stability transfer theorem in this setting, in a similar way as it was done in AECs (\cite{BaKuVa}). 

In this paper, we prove that if $\kappa$ is a strongly compact cardinal, then any MAEC with L\"owenheim-Skolem number below $\kappa$ is $<\kappa$-d-tame using ideas from \cite{Bo14} and \cite{Bo1x}.  Building on improvements from Brooke-Taylor and Rosick\'{y} \cite{BTRo} and Boney and Unger \cite{BoUn}, we actually weaken the hypothesis to an {\it almost strongly compact} (see Definition \ref{asc-def}); note that there is not a good notion of continuous infinitary logic\footnote{There are several approaches to defining this notion (see Eagle \cite{Eag} or Ben Yaacov and Iovino \cite{BYIo} for different frameworks), but neither seem to exactly capture the logic necessary for, say, a MAEC version of Shelah's Presentation Theorem.},so there is no loss here. After seeing early reports on this work, Lieberman and Rosick\'{y} informed us that they have achieved a similar result based on the ideas in \cite{LieRos}.

\section{Some basics on MAECs}

We consider a natural adaptation of the notion of {\it Abstract
Elementary Class} (see \cite{Gr} and \cite{BaMon}), but work in
a context of continuous logic that generalizes the continuous first-order setting of~\cite{CoMon} by removing the assumption of
uniform continuity\footnote{Uniform continuity guarantees logical
  compactness in their formalization, but we drop compactness in
  AEC-like settings.}. We base our definitions on \cite{HiHy,Za12}.
  
\begin{definition}\label{DensityCharacter}
The \emph{density character} of a topological space is the smallest
cardinality of a dense subset of the space. If $X$ is a topological
space, we denote its density character by $dc(X)$. If $A$ is a
subset of a topological space $X$, we define
$dc(A):=dc(\overline{A})$.
\end{definition}

\begin{definition}\label{MAEC}
Let $\mathcal{K}$ be a class of $L$-structures (in the context of
continuous logic) and $\ordencea$ be a binary relation defined in
$\C{K}$. We say that $(\C{K},\ordencea)$ is a {\it Metric Abstract
Elementary Class} (shortly {\it MAEC}) if:

\begin{enumerate}
\item $\C{K}$ and $\ordencea$ are closed under isomorphism.
\item $\ordencea$ is a partial order in $\C{K}$.
\item If $M\ordencea N$ then $M\subseteq N$.
\item (\emph{Completion of Union of Chains}) If $(M_i:i<\lambda)$ is a $\ordencea$-increasing chain then
    \begin{enumerate}
    \item the function symbols in $L$ can be uniquely interpreted on the completion of
        $\bigcup_{i<\lambda} M_i$ in such a way that
        $\overline{\bigcup_{i<\lambda}M_i} \in \C{K}$
    \item for each $j < \lambda$ , $M_j \ordencea
        \overline{\bigcup_{i<\lambda} M_i}$
    \item if each $M_i\ordencea  N$, then
        $\overline{\bigcup_{i<\lambda} M_i} \ordencea N$.
    \end{enumerate}
    \item (\emph{Coherence}) if $M_1\subseteq M_2\ordencea M_3$ and
        $M_1\ordencea M_3$, then $M_1\ordencea M_2$.
    \item (DLS) There exists a cardinality $LS(K)$ (which is called the
      {\it metric L\"o\-wen\-heim-Skolem number}) such that if $M
      \in \C{K}$
        and $A \subseteq M$, then there exists $N\in \C{K}$ such that $dc(N) \le dc(A) + LS(K)$
        and $A\subseteq N \ordencea M$.
\end{enumerate}
\end{definition}

Many basic notions in AECs ($\C{K}$-embeddings, amalgamation, etc.) translate directly to MAECs, so we omit their definition here; consult \cite{BaMon} or other sources.  Additionally, although not necessary, we work in the context of a monster model for simplicity; see \cite[Theorem 5.4]{Bo14} for how to transfer the concepts.





\begin{remark}[Monster Model]\label{Monster_Model}
If $\C{K}$ is an MAEC which satisfies AP and JEP and has large enough models, then we can
construct a very large model $\monster$ (which we call a {\it
  monster model}) which is homogeneous --i.e., every isomorphism
between two $\C{K}$-substructures of $\monster$ can be extended to an
automorphism of $\monster$-- and also universal --i.e., every model
with density character $<dc(\monster)$ can be $\C{K}$-embedded into
$\monster$.
\end{remark}

\begin{definition}[Galois type]\label{Galois_Type}
Under the existence of a monster model $\monster$ as in
Remark~\ref{Monster_Model}, for all $\tuple{a}\in \monster$ and
$N\ordencea \monster$, we define $\gatp(\tuple{a}/N)$ (the
\emph{Galois type of $\tuple{a}$ over $N$}) as the orbit of	
$\tuple{a}$ under $Aut(\monster/N):=\{f\in Aut(\monster): f\rest N =
id_{N}\}$. We denote the space of Galois types over a model $M\in
\C{K}$ by $\gaS(M)$.
\end{definition}

\begin{fact}
Assume there exists a monster model $\monster$.
Let $\C{M}\in \C{K}$ and $a,b\in \monster$. $\gatp(a/\C{M})=\gatp(b/\C{M})$ iff there exist $\C{M}_1,\C{M}_2\ordenceac \C{M}$, $\C{N}\in \C{K}$ and  $\C{K}$-embeddings $f_1:\C{M}_1\to \C{N}$ and $f_2:\C{M}_2\to \C{N}$ such that $a\in \C{M}_1$, $b\in \C{M}_2$, $f_1(a)=f_2(b) $ and $f_1\rest M = f_2\rest M =id_M$.
\end{fact}

The metric on the elements can be extended to a distance on types in a natural way.

\begin{definition}[Distance between types]\label{Distance_Types}
Let $p,q\in \gaS(M)$. We define $d(p,q):=\inf\{d(\tuple{a},\tuple{b}):
\tuple{a},\tuple{b}\in \monster, \tuple{a}\models p , \tuple{b}\models
q\}$, where ${\rm lg}(\tuple{a})={\rm lg}(\tuple{b})=:n$ and
$d(\tuple{a},\tuple{b}):=\max\{d(a_i,b_i): 1\le i\le n\}$.
\end{definition}

As defined, the above distance is just a pseudometric.  To turn it into a metric, the following assumption is used.

\begin{definition}[Continuity of Types]\label{Continuity_Types}
Let $\C{K}$ be an MAEC and consider $(a_n)\to a$ in $\monster$. We say that $\C{K}$ satisfies {\it Continuity of Types Property}\footnote{This property is also called {\it Perturbation Property in \cite{HiHy}}} (for short, {\it CTP}), if and only if, if $\gatp(a_n/M)=\gatp(a_0/M)$ for all $n<\omega$ then $\gatp(a/M)=\gatp(a_0/M)$.
\end{definition}

While typically used in work on MAECs, this paper does not require the assumption of CTP.

The following are the notions of tameness and type-shortness from AECs properly generalized to the continuous context of MAECs.  The first is from Zambrano \cite{Za12} and the second is new to this paper.

\begin{definition}[$\mathbf{d}$-tameness]\label{tameness}
Let $\C{K}$ be a MAEC and $\mu\ge LS(\C{K})$. We say that $\C{K}$ is {\it $<\mu$-$\mathbf{d}$-tame}
iff for every $\E>0,$ there exists $\delta_\E>0$ such that if for any $M\in \C{K}$ of density character $\ge\mu$ we have that $\dist(p,q)\ge \E$ where $p,q\in \gaS(M)$, then there exists $N\ordencea M$ of density character $<\mu$ such that $\dist(p\rest N,q\rest N)\ge \delta_\E$.
\end{definition}

\begin{definition}[$\mathbf{d}$-type shortness]\label{tameness}
Let $\C{K}$ be a MAEC and $\mu\ge LS(\C{K})$. We say that $\C{K}$ is {\it $\mu$-$\mathbf{d}$-type short}
iff for every $\E>0$ there exists $\delta_\E>0$ such that if for any $I$ of size $\geq \mu$ and $M\in \C{K}$, we have that $\dist(p,q)\ge \E$ where $p,q\in \gaS^I(M)$, then there exists $I_0 \subset I$ of size $\mu$ such that $\dist(p^{I_0},q^{I_0})\ge \delta_\E$.
\end{definition}

\section{Ultraproducts in MAECs}

The main ingredient in the proof is the use of very complete ultraproducts.  We use almost strongly compact cardinals, first isolated by Bagaria and Magidor \cite{BaMa}.

\begin{definition}[Almost strongly compact cardinal] \label{asc-def}
A cardinal $\kappa>\aleph_0$ is said to be {\it almost strongly compact} iff for any $\delta < \kappa$, any $\kappa$-complete filter can be extended to a $\delta$-complete ultrafilter. 
\end{definition}

\begin{fact}\label{StronglyCompact}
The following are equivalent:
\begin{enumerate}
\item $\kappa$ is almost strongly compact.
\item For any $\delta < \kappa$, both (discrete) $L_{\delta,\omega}$ and $L_{\delta, \delta}$ satisfy $\kappa$-compactness.
\item For every $\delta < \kappa \leq \lambda$, there is a fine, $\delta$-complete ultrafilter $U$ on $P_{\kappa}\lambda$ (i.e.,  a $\delta$-complete ultrafilter such that for every $\alpha<\kappa$, we have $[\alpha] = \{X \in P_{\kappa}\lambda : \alpha \in X\} \in U$).
\end{enumerate}
\end{fact}

We want to have an analogue of \L o\'{s}' Theorem for AECs in the metric setting.  Rather than reproving a continuous version of \L o\'{s}' Theorem, we make use of the fact that there is an ``invertible functor'' from the class of Metric AECs to the class of AECs.

\begin{theorem}[Boney, Theorem 6.1 of \cite{Bo1x}]\label{Kdense}
Let $L$ be a continuous language. Then there is a discrete language $L^+$ such that, for every MAEC $\C{K}$ with $LS(\C{K}) =|L|$, there is
\begin{enumerate}
\item an AEC $\C{K}_{dense}$ with $L(\C{K}_{dense}) =|L^+|$ and $LS (\C{K}_{dense}) =LS (K)$;
\item  a map from $M\in \C{K}$ and nicely dense subsets $A$ of $M$ to $M_A \in \C{K}_{dense}$; and
\item  a map from $A\in \C{K}_{dense}$ to $\overline{A}\in \C{K}$
\end{enumerate}
with the properties that
\begin{enumerate}
\item $M_A$ has universe $A$ 
\item $\overline{A}$ has universe that is the completion of $A$ with respect to the derived metric 
\item  The maps above are essentially inverses (up to the choice of limit)
\item
	\begin{enumerate}
	\item Given $M_l \in \C{K}$ and $A_l$ nicely dense in $M_l$ for $l= 0,1$, if $f:M_0  \to M_1$ is a $\C{K}$ -embedding such that
	$f[A_0] \subset A_1$, then $f\upharpoonright A_0$ is a $K_{dense}$-embedding from $(M_0)_{A_0}$ to $(M_1)_{A_1}$.
	\item Given $A,B\in\C{K}_{dense}$ and a $K_{dense}$-embedding $f:A\to B$, this lifts canonically to a $K$-embedding $f:\overline{A}\to \overline{B}$
	\end{enumerate}
\end{enumerate}
\end{theorem}

For the sake of completeness, we provide the definition of ultraproduct of complete metric structures.  In using countably complete ultrafilters, we have an advantage over the normal situation.  Typically, one considers only the ultraproduct of uniformly bounded metric spaces or considers them as pointed metric spaces (and the choice of point can affect the ultraproduct).  However, with a countably complete metric space, none of this is necessary.

\begin{definition}
Let $X$ be a topological space and let $(x_i)_{i\in I}$ be a family of elements of $X$. If $U$ is an ultrafilter on $I$ and $x\in X$, 
$lim_{U}x_i =x $ means that for every neighborhood $V$ of $x$, $\{i \in I : x_i \in  V\}\in U$.
\end{definition}

\begin{definition}
Let $\langle (M_i,d_i) : i \in I\rangle $ be a family of metric spaces. Let $U$ be a countably complete ultrafilter on $I$. Define $d$ on $\prod_{i\in I} M_i$ by $d((x_i), (y_i)): = lim_U d_i(x_i, y_i)$. Define $(x_i)\sim_U (y_i)$ iff $d((x_i),(y_i))=0$.
\end{definition}

\begin{fact}
$(\prod_{i\in I} M_i)/\sim_U,d)$ is a complete metric space.
\end{fact}

The proof of this fact is mostly standard (see \cite{CoMon}).  The nonstandard part is the use of countable completeness to remove the uniformly bounded assumption.  Let $(x_i)_U, (y_i)_U \in \prod_{i \in I} M_i/\sim_U$; we want to show that $d((x_i)_U, (y_i)_U)$ is a real number.  It is enough to show that there is some $N$ such that $\{i \in I : d_i(x_i, y_i) < N\} \in U$.  We can partition $I$ into countably many pieces as $X_N := \{ i \in I : N \leq d_i(x_i, y_i) < N+1\}$.  Since $U$ is countably complete, some $X_N \in U$ and this gives us the bound.

\begin{fact}
Let $\langle \mathcal{M}_i : i \in I\rangle $ be a sequence of $L$-structures in continuous logic.  $(\prod_{i\in I} M_i)/\sim_U,d)$ together with all the ultraproduct of  interpretations of symbols in $L$ is an $L$-structure in Continuous Logic (denoted by $\prod^d_U {M}_i$).
\end{fact}

In the following lines, we prove the \L o\'{s}' theorem for MAECs, as in \cite{Bo14}.

\begin{theorem}[\L o\'{s}' Theorem for MAECs (1)]\label{Los1}
If $\langle M_i :i\in I\rangle$ is a sequence of structures in $\C{K}$ and $U$ is an $LS(K)^+$-complete ultrafilter, then the metric ultraproduct $\prod^d_{U} M_i$ belongs to $\C{K}$.
\end{theorem}

In order to prove this, we need to see that the metric and discrete ultraproducts agree given enough completeness.

\begin{lemma}\label{MDultraproduct}
Suppose that $A_i$ is nicely dense in $M_i \in K$ for each $i \in I$ and $U$ is a countably complete ultraproduct on $I$.  Then
$$\overline{\prod_U M_{A_i}} = \prod_U^d M_i$$
where $\prod^d$ denotes the ultraproduct as a metric structure.
\end{lemma}

\bdem

The key part of the proof is the above argument that countable completeness simplifies the metric ultraproduct of unbounded structures.  
We show that the universes of the two models are the same (up to canonical isomorphism); that they have the same structure follows similarly.

Let $(x_i)_U \in \prod^d_U M_i$.  For each $n < \omega$, there is $x_i^n \in A_i$ such that $d_i(x_i^n, x_i) < \frac{1}{2^n}$.  Then set $x^n := (x_i^n)_U \in \prod_U M_{A_i}$.  Then, for $n < m$, we have
$$\left\{ i : M_{A_i} \vDash d(x_i^n, x_i^m) < \frac{1}{2^{n-1}} \right\} \in U$$
Thus, $(x^n)_n$ is a Cauchy sequence in $\prod_U M_{A_i}$, so it has a limit in $x_* \in \overline{\prod_U M_{A_i}}$.  It is clear that this limit must be $(x_i)_U$ (up to canonical isomorphism fixing $\prod_U M_{A_i}$).

Now suppose $x \in \overline{\prod_U M_{A_i}}$.  This means that there is $(x^n_i)_U \in \prod_U M_{A_i}$ such that, for each $n < \omega$, $d_{\prod_U M_{A_i}}( x, (x_i^n)_U) < \frac{1}{2^n}$.  By the countable completeness of $U$, there is $X \in U$ such that, for all $i \in X$, the sequence $\{ (x_i^n)_U: n < \omega\}$ is Cauchy.  Since $M_i$ is complete, this converges to $y_i \in M_i$.  Then $(y_i)_U$ (which is defined on a $U$-large set) is the limit of $\left((x_i^n)_U\right)_n$.  Thus, $x = (y_i)_U$ up to a canonical isomorphism.
\edem
\bdem[Proof Thm.~\ref{Los1}]

Let $A_i\subset M_i$ ($i\in I$) be nicely dense subsets.  By Theorem \ref{Kdense}, $M^i:=(M_i)_{A_i} \in \C{K}_{dense}$.  Note that $LS(\C{K}_{dense}) = LS(\C{K})$.  Thus, by \L o\'{s}' Theorem for discrete AECs \cite[Theorem 4.3]{Bo14}, we have that $\prod_U M^i \in \C{K}$.  Thus, by Theorem \ref{Kdense} again, $\overline{\prod_U M^i} \in \C{K}$.  By Lemma \ref{MDultraproduct}, we have $\prod_U^d M_i \in \C{K}$ as desired.

\edem



We can transfer other clauses from \cite[Theorem 4.3]{Bo14} similarly:

\begin{theorem}[\L o\'{s}' theorem for MAECs (2)]
\label{Los2}\label{UltraproductIsos}
Let $U$ be an $LS(\C{K})^+$-complete ultrafilter over $I$.
\begin{enumerate}

	\item If $\langle M_i : i\in I\rangle$ and $\langle N_i : i\in I\rangle$ are sequences of structures in a MAEC $\C{K}$ such that $M_i\ordencea N_i$ for every $i\in I$, then $\prod_U^d M_i \ordencea \prod_U^d N_i$.
	
	\item Let $\langle M_i:i\in I\rangle$ and $\langle N_i:i\in I\rangle$ be sequences of structures in $\C{K}$ and $\{h_i : i\in I\}$ a family of $L(\C{K})$-isomorphisms $h_i:M_i\stackrel{\cong}{\to} N_i$, then the mapping $h: \prod_U^d M_i \to \prod_U^d N_i$ defined by $h((x_i)_U):= (h_i(x_i))_U$ is an $L(\C{K})$-isomorphism. 

	\item Let $\langle M_i : i\in I\rangle$ and $\langle N_i : i\in I\rangle$ be sequences of structures in $\C{K}$ and $\{h_i:M_i\to N_i:i\in I\}$ be a family of $\C{K}$-embeddings. Then $\prod_U h_i : \prod_U^d M_i\to \prod_U^d N_i$ defined in theorem~\ref{UltraproductIsos} is a $\C{K}$-embedding.

\end{enumerate}	
\end{theorem}

Now we have all the tools to prove our main theorem.

\begin{theorem}
If $\C{K}$ is a MAEC with $LS(\C{K}) < \kappa$ with $\kappa$ almost strongly compact, then $K$ is $<\kappa-d$-tame and $<\kappa-d$-type short.
\end{theorem}


\bdem
Let $\C{K}$ be an MAEC with $LS(\C{K}) < \kappa$.  We want to show that $\C{K}$ is $<\kappa-d$-tame; the proof of $<\kappa-d$-type shortness is similar. We will show the contrapositive of the statement from the definition of tameness.

Let $\varepsilon>0$ and $0<\delta_\varepsilon \leq \varepsilon$. Let $p, q \in S(M)$ such that, for all $M^- \in P^*_\kappa M$, $d(p \rest M^-, q \rest M^-) < \delta_\varepsilon$.  Let $a \vDash p$ and $b \vDash q$
.  For each $M^- \in P^*_\kappa M$, find the following:
\begin{itemize}
	\item $a_{M^-} \vDash p \rest M^-$ and $b_{M^-} \vDash q \rest M^-$ such that $d(a_{M^-}, b_{M^-}) < \delta_\varepsilon 
	$; and
	\item $f_{M^-}, g_{M^-} \in Aut_{M^-} \sea$ such that $f_{M^-}(a_{M^-}) = a$ and $g_{M^-}(b_{M^-}) = b$.
\end{itemize}
By almost strong compactness, find a $LS(K)^+$-complete, fine ultrafilter on $P^*_\kappa M$.  Now, take the ultraproduct by the \L o\'{s}' Theorem and get:
\begin{itemize}
	\item the average $f \in Aut \prod_U \sea$ takes $a^*:=(M^- \mapsto a_{M^-})_U$ to $ (M^- \mapsto a)_U$;
	\item the average $g \in Aut \prod_U \sea$ takes $b^*:=(M^- \mapsto b_{M^-})_U$ to $(M^- \mapsto b)_U$;
	\item the ultrapower embedding $h$ takes $a$ to $(M^- \mapsto a)_U$ and $b$ to $(M^- \mapsto b)_U$;
	\item $d(a^*, b^*) = \lim_U d(a_{M^-}, b_{M^-}) < 
	\delta_\varepsilon $ (we get ``$\leq \delta_\varepsilon$'' from general facts about $U$-limits and strengthen it to strict inequality by countable completeness); and
	\item $f$ and $g$ fix $h(M)$ by fineness as in \cite{Bo14}: let $m \in M$.  Then $\{M^- \in P^*_\kappa M : m \in M^-\} \in U$ and, for each such $M^-$, $f_{M^-} (m) = m$.  Thus, $f(h(m)) = (M^- \mapsto f_{M^-}(m))_U = (M^- \mapsto m)_U = h(m)$.  Similarly for $g$.
\end{itemize}
Thus, we have that $a^* \vDash h(p)$ and $b^* \vDash h(q)$ and $d(a^*, b^*) < \delta_\varepsilon
$.  This gives $d(h(p), h(q)) <
\delta_\varepsilon 
\leq\varepsilon$.  Since $h$ is a $\C{K}$-embedding
, this means that $d(p, q) < \epsilon$.

\edem

As mentioned above, this result holds even without a monster model.  One important caveat is that, without even amalgamation, the result holds about atomic equivalence of triples (i. e., before the transitive closure is taken).  However, when the transitive closure is taken, the result might no longer paper.

Further more, one can extend the other results of \cite{Bo14} about measurable cardinals and weakly compacts to this setting.

\end{document}